\documentclass[12pt]{article}
\usepackage{amsmath, amsthm, amssymb, amstext}
\usepackage{mathrsfs}
\usepackage{array, amsfonts, mathrsfs}
\usepackage{latexsym, amsmath}
\usepackage{graphicx,color}
\usepackage{epsf}

\newtheorem{Theorem}{Theorem}

\date{}

\begin{document}

\author{M.I.Belishev\thanks {St. Petersburg Department of Steklov Mathematical
        Institute, St.Petersburg, Russia, e-mail: belishev@pdmi.ras.ru.},\,
        D.V.Korikov\thanks {St.Petersburg Department of Steklov Mathematical
        Institute, St. Petersburg, Russia, e-mail:
        thecakeisalie@list.ru.}}

\title{On characterization of Dirichlet-to-Neumann map of Riemannian surface with boundary}
\maketitle

\begin{abstract}
Let $(M,g)$ be a smooth compact orientable two-dimensional
Riemannian manifold ({\it surface}) with a smooth metric tensor
$g$ and smooth connected boundary $\Gamma$. Its {\it DN-map}
$\Lambda_g:{C^\infty}(\Gamma)\to{C^\infty}(\Gamma)$ is associated
with the (forward) elliptic problem $ \Delta_gu=0 \,\,\, {\rm
in}\,\,M\setminus\Gamma,\,\,u=f \,\,\, {\rm on}\,\,\,\Gamma$, and
acts by
 $
\Lambda_g f:=\partial_\nu u^f \,\,\, {\rm on}\,\,\,\Gamma,
 $
where $\Delta_g$ is the Beltrami-Laplace operator, $u=u^f(x)$ is
the solution, $\nu$ is the outward normal to $\Gamma$. The
corresponding {\it inverse problem} is to determine the surface
$(M,g)$ from its DN-map $\Lambda_g$.

We provide the necessary and sufficient conditions on an operator
acting in ${C^\infty}(\Gamma)$ to be the DN-map of a surface. In
contrast to the known conditions by G.Henkin and V.Michel in terms
of multidimensional complex analysis, our ones are based on the
connections of the inverse problem with commutative Banach
algebras.
\end{abstract}

\subsubsection*{0.\,\,Introduction}
\noindent$\bullet$\,\,\,Let $(M,g)$ be a smooth
\footnote{Throughout the paper {\it smooth} means
$C^\infty$-smooth} compact orientable two-dimensional Riemannian
manifold with a smooth metric tensor $g$ and smooth connected
boundary $\Gamma$. In what follows, we deal with the manifolds of
this class only and, for short, call them the {\it surfaces}. The
Dirichlet-to-Neumann operator ({\it DN-map})
$\Lambda:{C^\infty}(\Gamma)\to{C^\infty}(\Gamma)$ of the surface
is associated with the (forward) elliptic problem
 \begin{align}
\label{Eq 1} & \Delta_gu=0 \qquad {\rm in}\,\,M\setminus\Gamma, \\
\label{Eq 2} & u=f \qquad\quad\, {\rm on}\,\,\,\Gamma
 \end{align}
and acts by the rule
 $$
\Lambda_g f:=\partial_\nu u^f \qquad {\rm on}\,\,\,\Gamma,
 $$
where $\Delta_g$ is the Beltrami-Laplace operator, $u=u^f(x)$ the
solution, $\nu$ the outward normal to $\Gamma$. The corresponding
{\it inverse problem} is to recover the surface $(M,g)$ via the
operator $\Lambda_g$. In applications it is also known as the
Electric Impedance Tomography problem. More generally, one needs
to answer the question: to what extent does the DN-map determine
the surface?

In the paper \cite{LUEns} by M.Lassas and G.Uhlman, it is shown
that the DN-map determines the surface $M$ up to conformal
equivalence. In more detail, if $(M,g)$ and $(M',g')$ have the
common boundary $\Gamma$, and $\Lambda_g=\Lambda_{g'}$ holds,
   then there exists a diffeomorphism $\psi: M'\to M$ provided
$\psi|_\Gamma=\rm id$ and a smooth positive function $\rho$
obeying $\rho|_\Gamma\equiv 1$, such that $g=\rho\,\psi_*g'$.

In the paper \cite{BCald} by M.I.Belishev, the same result is
established by the use of the connections between the EIT problem
and the holomorphic function algebra of the surface. Moreover, the
formulas, which express the topological invariants of the surface
(Betti numbers) in terms of the DN-map, are provided. In
\cite{BSharaf} these formulas are generalized on the
multidimensional case. The paper \cite{BKor_JIIPP} extends the
algebraic approach to nonorientable surfaces.
\smallskip

\noindent$\bullet$\,\,\,All the above results refer to the
situation when the operator $\Lambda_g$ {\it is given}, and the
existence of the surface, for which $\Lambda_g$ is the DN-map, is
a priori assumed. Thus, for such a $\Lambda_g$, the solvability of
the EIT problem is guaranteed. However, an important question
remains about the {\it inverse data characterization}, i.e., on
the necessary and sufficient conditions for an operator $\Lambda$
to be the DN-map of a surface (to satisfy $\Lambda=\Lambda_g$). In
other words, we are talking about a criterion for solvability of
the EIT problem.

Such a criterion is presented in the paper \cite{H&M} by
G.M.Henkin and V.Michel in terms of multidimensional complex
analysis. In the given paper, we propose a characterization based
on the connections of EIT problem with Banach algebras. So, the
novelty is a new formulation of the solvability conditions. The
list of our conditions is also rather long  but, however, we would
venture to claim that our formulation is more transparent and
simpler than that proposed in \cite{H&M}.

Our approach makes the use of the classical result \cite{Aupetit}
on the existence of a complex structure on the Gelfand spectrum of
a commutative Banach algebra. It is the result, which provides the
sufficiency of the proposed characteristic conditions.
\smallskip

\noindent$\bullet$\,\,\, In the first section, the list of the
necessary and sufficient conditions for solvability of the EIT
problem is specified. The second section contains the proof of
necessity of these conditions, the proof being based on some
general properties of the DN-map. In the third section, the
sufficiency is proved. For the convenience of the reader, we
present the basic definitions and minimal required information on
the Banach algebras.

\subsubsection*{1.\,\,Main result}

\noindent$\bullet$\,\,\,Let $\Gamma$ be a smooth curve
diffeomorphic to a circle, $d\gamma$ its length element, $\gamma$
a contininuous tangent field of unit vectors on $\Gamma$, and
$\Lambda: \ C^{\infty}(\Gamma;\mathbb{R})\mapsto
C^{\infty}(\Gamma;\mathbb{R})$ a linear map. With $\Lambda$ one
associates the map $\Upsilon : \
C^{\infty}(\Gamma;\mathbb{C})\mapsto
C^{\infty}(\Gamma;\mathbb{C})$,
 \begin{equation}
\label{Goperator}
\Upsilon\zeta:=(\Lambda\Re\zeta-\partial_{\gamma}\Im\zeta)+i(\Lambda\Im\zeta+\partial_{\gamma}\Re\zeta).
 \end{equation}
As is easy to verify, it is a (complex) linear operator. Also, for
$\eta\in C^{\infty}(\Gamma;\mathbb{C})$ and
$z\in\mathbb{C}\backslash\eta(\Gamma)$, introduce the map
$\Upsilon_{\eta,z} : \ C^{\infty}(\Gamma;\mathbb{C})\mapsto
C^{\infty}(\Gamma;\mathbb{C})$ as
\begin{equation}
\Upsilon_{\eta,z}\zeta:=\Upsilon\frac{\zeta}{\eta-ze}\,,
\end{equation}
where $e$ is the function equal to $1$ on $\Gamma$.

Let $I$ be the identity operator on
$C^{\infty}(\Gamma;\mathbb{R})$,
$\partial_{\gamma}C^{\infty}(\Gamma;\mathbb{R})$ the space of
smooth real-valued functions with zero mean value on $\Gamma$, $J:
\ \partial_{\gamma}C^{\infty}(\Gamma;\mathbb{R})\mapsto
\partial_{\gamma}C^{\infty}(\Gamma;\mathbb{R})$ the integration on
$\Gamma$: $J\partial_{\gamma}=\partial_{\gamma}J=I$. By $\sharp S$
we denote the cardinality of $S$.

\smallskip

Our main result is the following.
 \begin{Theorem}\label{main}
The operator $\Lambda$ is the DN-map of a surface if and only if
it satisfies the conditions:
\smallskip

\noindent{\bf i.}\,\,\,$e\in{\rm Ker}\Upsilon$ and
$\zeta_{1}\zeta_{2}\in {\rm Ker}\Upsilon$ for any
$\zeta_{1},\zeta_{2}\in {\rm Ker}\Upsilon$;
\smallskip

\noindent{\bf ii.}\,\,\,if $\zeta_{1},\zeta_{2}\in {\rm
Ker}\Upsilon$, $\zeta_{1}/\zeta_{2}\in
C^{\infty}(\Gamma;\mathbb{C})$, and there exists a polynomial $P$,
${\rm deg\,}P\geqslant 1$ such that $P(\zeta_{1}/\zeta_{2})\in
{\rm Ker}\Upsilon$, then $\zeta_{1}/\zeta_{2}\in {\rm
Ker}\Upsilon$;
\smallskip

\noindent{\bf iii.}\,\,\,$\overline{{\rm Ker}\Upsilon}\cap
C^{\infty}(\Gamma;\mathbb{C})={\rm Ker}\Upsilon$\,\,\,(the closure
in $C(\Gamma;\mathbb{C})$);
\smallskip

\noindent{\bf iv.}\,\,\,${\rm dim}(\partial_{\gamma}+\Lambda
J\Lambda)C^{\infty}(\Gamma;\mathbb{R})<\infty$;
\smallskip

\noindent{\bf v.}\,\,\,if $\eta\in {\rm Ker}\Upsilon$ and
$z\in\mathbb{C}\backslash\eta(\Gamma)$, then
\begin{equation}\label{algebraic argument principle}
{\rm dim\,}[\Upsilon_{\eta,z}{\rm Ker}\Upsilon]\,=\,\frac{1}{2\pi
i}\int\limits_{\Gamma}\frac{\partial_{\gamma}\eta}{\eta-ze}d\gamma\,;
\end{equation}
\smallskip

\noindent{\bf vi.}\,\,\,for any $x\in\Gamma$, there exist a
function $\eta_{x}\in {\rm Ker}\Upsilon$ and a neighborhood
$U_{x}\ni\eta_{x}(x)$ diffeomorphic to an open disk
$D\subset\mathbb C$, such that
\begin{enumerate}
\item $\partial_{\gamma}\eta_{x}(x)\ne 0$ is valid and there is no
points on $\Gamma$, at which all derivatives
$\partial_{\gamma}^{k}\eta_{x}$, $k\geqslant 1$ vanish
simultaneously, whereas $\sharp\,\eta_{x}^{-1}(\{z\})<\infty$
holds for all $z\in\mathbb C$;

\item $\Upsilon_{\eta_{x},z}e=0$ holds on one connected
component of $U_{x}\backslash\eta_x(\Gamma)$, whereas
$\Upsilon_{\eta_{x},z}e\ne 0$ holds on the other connected
component; \item the equation
\begin{equation}
\label{boundeq} \Upsilon\Big(\frac{\zeta-ce}{\eta_{x}-ze}\Big)=0
\text{ on } \Gamma
\end{equation}
has a solution $c\in\mathbb{C}$ for any $z\in U_{x}$ and
$\zeta\in{\rm Ker}\Upsilon$;
\end{enumerate}
\smallskip

\noindent{\bf vii.}\,\,\,if $\zeta,1/\zeta\in{\rm Ker}\Upsilon$,
then $\Lambda\,{\rm log\,}|\zeta|=\partial_{\gamma}{\rm
arg\,}\zeta$.
 \end{Theorem}
As a comment, note the following. Condition $\bf i$ means that
${\rm Ker}\Upsilon$ is an algebra, whereas\, $\bf vi$\, shows that
this algebra must be rich enough to contain the functions $\eta_x$
with the required properties. By condition ${\bf vi}.2$, the
function $\eta_{x}-ze$ is or is not invertible in the algebra
${\rm Ker}\Upsilon$ depending on the position of $z$ on the
complex plane. Condition ${\bf vi}.3$, from an algebraic point of
view, means that $(\eta_x-ze)\,{\rm Ker}\Upsilon$ is an ideal in
${\rm Ker}\Upsilon$ of codimension 1, i.e., is a maximal ideal.
Also, as is easy to see, the embedding $g\in{\rm Ker}\Upsilon$
implies $\partial_{\gamma}\Re g,\partial_\gamma\Im g\in{\rm
Ker\,}[I+(\Lambda J)^2]$. The operator $I+(\Lambda J)^2$ is a key
object of the papers \cite{BCald,BKor_JIIPP,BSharaf}.
\medskip

The rest of the paper is devoted to the proof of Theorem
\ref{main}.

\subsubsection*{2.\,\,Necessity}
Here we show that any DN-map satisfies conditions {\bf i}--{\bf
vii}.
\smallskip

\noindent$\bullet$\,\,\,
Suppose that $\Lambda=\Lambda_g$ is the DN-map of some surface
$(M,g)$. Recall that $\gamma$ and $\nu$  are the tangent and
normal unit vector fields at the boundary $\Gamma$.

Choose a continuous family of rotations $M\ni x\mapsto\Phi_{x}\in
{\rm End\,}T_{x}M$,
$$g(\Phi_{x}a,\Phi_{x}b)=g(a,b), \quad g(\Phi_{x} a,a)=0, \qquad a,b\in T_{x}M, \ x\in M,$$
such that $\Phi\nu=\gamma$ on $\Gamma$. A function $w\in
C^{\infty}(M;\mathbb{C})$ is called {\it holomorphic} if the
Cauchy-Riemann condition $\nabla_g\Im w=\Phi\nabla_g\Re w$ holds
in $M$. Let $w$ be holomorphic and $\zeta=w|_\Gamma$ be its trace
on the boundary. The real functions $\Re w$ and $\Im w$ are
harmonic in ${\rm int}M$ and provide the solutions $\Re
w=u^{\Re\zeta}$ and $\Im w=u^{\Im\zeta}$ to (\ref{Eq 1}), (\ref{Eq
2}). Restricting the Cauchy-Riemann conditions on $\Gamma$, one
obtains
\begin{equation}
\label{CRonboundary1} \Lambda\Re\zeta=\partial_{\nu}\Re
w=\partial_{\gamma}\Im\zeta , \qquad
\Lambda\Im\zeta=\partial_{\nu}\Im w=-\partial_{\gamma}\Re\zeta
\quad \text{ on } \Gamma,
\end{equation}
which implies $\Upsilon(\zeta)=0$ according to (\ref{Goperator}).

Now, suppose that $\zeta\in C^{\infty}(\Gamma;\mathbb{C})$ and
$\Upsilon(\zeta)=0$. Then the function
$w:=u^{\Re\zeta}+iu^{\Im\zeta}$ is holomorphic in ${\rm int}M$.
Indeed, since $\Upsilon(\zeta)=0$, one has (\ref{CRonboundary1}),
i.e., $\nabla_g\Im w=\Phi\nabla_g\Re w$ holds on $\Gamma$. Let $U$
be an arbitrary neighborhood in $M$ diffeomorphic to the disc, and
$\partial U\cap\Gamma$ contains a segment $\Gamma'$ of non-zero
length. Since $\partial_{\nu}\Re w=u^{\Re\zeta}$ is harmonic in
$U$, there exists a function $v$ such that $\nabla_g
v=\Phi\nabla_g\Re w$ in $U$. Thus, $\partial_{\nu}\Re
w=\partial_{\gamma}v$ and $\partial_{\nu}v=-\partial_{\gamma}\Re
w$ on $\Gamma'$. Comparing with (\ref{CRonboundary1}), one obtains
$v=\Im w+{\rm const}$, $\partial_{\nu}v=\partial_{\nu}\Im w$ on
$\Gamma'$. So, $v$ and $\Im w+{\rm const}$ are harmonic in $U$ and
have the same Cauchy data on $\Gamma'$. Due to the uniqueness of
the solution to the Cauchy problem for the second order elliptic
equations, $v$ coincides with $\Im w+{\rm const}$ in $U$, and
$\nabla_g\Im w=\nabla_g v=\Phi\nabla_g\Re w$ in $U$. Since $U$ is
arbitrary, $\nabla_g\Im w=\Phi\nabla_g\Re w$ holds in $M$, and $w$
is holomorphic. So, we have proved that ${\rm Ker}\Upsilon$
coincides with the set of traces on $\Gamma$ of all holomorphic
smooth functions on $M$. Obviously, such a set is an algebra with
respect to point-wise multiplication and $e\in {\rm Ker}\Upsilon$.
So, {\bf i} is valid.
\smallskip

\noindent$\bullet$\,\,\,
Let $\zeta_{1},\zeta_{2}\in {\rm Ker}\Upsilon$,
$\zeta=\zeta_{1}/\zeta_{2}\in C^{\infty}(\Gamma;\mathbb{C})$, and
$P(\zeta)\in {\rm Ker}\Upsilon$, where $P$ is a polynomial of
degree $p>1$. In view of the already proven, there exist
holomorphic functions $w_{1},w_{2},w_{P}$ such that
$w_{1}|_{\Gamma}=\zeta_{1}$, $w_{2}|_{\Gamma}=\zeta_{2}$, and
$w_{P}|_{\Gamma}=P(\zeta)$. Then the function $w:=w_{1}/w_{2}$ is
meromorphic in ${\rm int}M$ and $w|_{\Gamma}=\zeta\in
C^{\infty}(\Gamma;\mathbb{C})$. The last implies that the poles of
$w$ do not accumulate to $\Gamma$ and the number of them is
finite. The function $P(w)$ is also meromorphic and its poles
coincide with those of $w$, while their multiplicities are $p$
times greater than those of $w$. Since $P(w)=P(\zeta)=w_{P}$ on
$\Gamma$, the function $P(w)$ coincides with $w_{P}$ outside the
poles of $w$ due to uniqueness of analytic continuation. Then
$P(w)=w_{P}$ everywhere on $M$. Thus, $w$ is holomorphic and its
trace $\zeta$ belongs to ${\rm Ker}\Upsilon$. This proves {\bf
ii}.
\smallskip

\noindent$\bullet$\,\,\,
Due to the maximum principle, the set $\overline{{\rm
Ker}\Upsilon}$ coincides with the set of traces on $\Gamma$ of all
holomorphic continuous functions. Since the smooth elements of
$\overline{{\rm Ker}\Upsilon}$ are the traces of holomorphic
smooth functions, one obtains {\bf iii}.

The property {\bf iv} follows from the equality
\begin{equation}
\label{dimhandles} {\rm dim}(\partial_{\gamma}+\Lambda
J\Lambda)C^{\infty}(\Gamma;\mathbb{R})=1-\mathscr{X}(M)
\end{equation}
(see formula (1.6), \cite{BCald}), where $\mathscr{X}(M)$ is the
Euler characteristics of $M$.
\smallskip

\noindent$\bullet$\,\,\,
Suppose that $\eta\in {\rm Ker}\Upsilon$ and
$z\in\mathbb{C}\backslash\eta(\Gamma)$. Then there exists the
holomorphic in ${\rm int}M$ function $w_{0}$ such that
$w_{0}|_{\Gamma}=\eta$. Denote by $x_{1},\dots,x_{l}$ all the
zeroes of $w_{0}-z$ and by $m_{1},\dots,m_{l}$ their
multiplicities. We make use of the argument principle
\footnote{for a compact Riemann surface with boundary, the
argument principle can be obtained by simple modification of the
proofs of Theorem 3.17 and Corollary 3.18, \cite{Miranda}.}:
\begin{equation}
\label{argprinc} \frac{1}{2\pi
i}\int\limits_{\Gamma}\frac{\partial_{\gamma}\eta}{\eta-ze}d\gamma=\sum_{k=1}^{l}m_{k}.
\end{equation}
Since $\Gamma$ is smooth, the manifold $(M,g)$ can be embedded
into a larger non-compact smooth manifold $(M',g')$, $g'|_{M}=g$.
For each $k=1,\dots,l$ and $s=0,\dots,m_{k}-1$, choose a
holomorphic into $M'$ function $w_{k,s}$ such that
$x_{1},\dots,x_{l}$ are all zeroes of $w_{k,s}$ on $M'$ and
multiplicity of $x_{j}$ is equal to $s$ if $j=k$ and to $m_{j}$ if
$j\ne k$. The existence of such $w_{k,s}$ follows from Proposition
26.5, \cite{Forster}. The linear combination
$\sum_{k,s}\frac{c_{k,s}w_{k,s}}{w_{0}-z}$ has no poles in $M$
only if all $c_{k,s}$ equal to zero. Denote
$\eta_{k,s}:=w_{k,s}|_{\Gamma}$; then
$\Upsilon_{\eta,z}(\sum_{k,s}c_{k,s}\eta_{k,s})= 0$ only if all
$c_{k,s}$ are zeros. Hence, the functions
\begin{equation}
\label{bas} \Upsilon_{\eta,z}(\eta_{k,s}), \quad k=1,\dots,l, \ \
s=0,\dots,m_{k}-1
\end{equation}
are linearly independent.

Now, suppose that $w\in C^{\infty}(\Gamma;\mathbb{C})$ is
holomorphic in $M$ and $\zeta=w|_{\Gamma}$. For any $k=1,\dots,
l$, there exist $d_{k,s}\in\mathbb{C}$ ($l=1,\dots,m_{k}$) such
that $w-\sum_{s}d_{k,s}w_{k,s}$ has a zero of multiplicity
$\geqslant m_{k}$ at $x=x_{k}$. Since $x_{k}$, $k'\ne k$ is a zero
of multiplicity $m_{k}$ for each each $w_{k',s}$, the function
$w-\sum_{k,s}d_{k,s}w_{k,s}$ have at each $x_{k}$ a zero of
multiplicity $\geqslant m_{k}$. Therefore, the ratio
$\frac{w-\sum_{k,s}d_{k,s}w_{k,s}}{w_{0}-z}$ is holomorphic in $M$
and, hence,
$\Upsilon_{\eta,z}(\zeta-\sum_{k,s}d_{k,s}\eta_{k,s})=0$. This
means that (\ref{bas}) is a basis in $\Upsilon_{\eta,z}{\rm
Ker}\Upsilon$. In particular, ${\rm dim}\Upsilon_{\eta,z}{\rm
Ker}\Upsilon=\sum_{k=1}^{l}m_{k}$. Comparing with
(\ref{argprinc}), one arrives at {\bf v}.
\smallskip

\noindent$\bullet$\,\,\,
Let $x$ be an arbitrary point of $\Gamma$. According to
Proposition 26.5, \cite{Forster}, there exists the holomorphic in
$M'$ function $w_{x}$ such that $x$ is unique zero of $w_{x}$ and
its multiplicity is equal to one. For any $c\in\mathbb{C}$, the
function $w_{x}-c$ has only finite number of zeros on $M$
(otherwise, there would be an accumulation point of such zeros due
to the compactness of $M$) and each zero of $w_{x}-c$ is of finite
multiplicity. This implies {\bf vi.}1 for the function
$\eta_{x}:=w_{x}|_{\Gamma}\in {\rm Ker}\Upsilon$.

Next, since $\nabla\Re w_{x}(x)\ne 0$, the map $w_{x}: M\to\mathbb
C$ is a bijection of a neighborhood $V_{0}$ of $x$ and
neighborhood $w_{x}(V_{0})$ of the zero, and $|w_{x}(x')|>0$ holds
for any $x'\in M'\backslash\{x\}$. Let $K$ be a compact in $M'$
that contains $M\cup V_{0}$. Then the set $K\backslash V_{0}$ is
also compact and $|w_{x}(x')|>c_{0}>0$ for any $x'\in K\backslash
V$. Choose a neighborhood $V_{1}\subset V_{0}$ sufficiently small
to obey $|w_{x}(x')|<c_{0}/2$ for any $x\in V_{1}$. Then the
pre-image $w_{x}^{-1}(\{z\})$ of any $z\in w_{x}(V_{1})$ is
contained in $V_{0}$ and, since $w_{x}$ is a bijection of $V_{0}$
and $w_{x}(V_{0})$, it consists of a single element. Denote
$U_{x}:=w_{x}(V_{1})$, $U_{x,1}:=U_{x}\backslash w_{x}(M)$, and
$U_{x,2}:=U_{x}\cap w_{x}(M)$. The function $\frac{1}{w_{x}-z}$
has no poles on $M$ for any $z\in U_{x,1}$ and has a simple pole
on $M$ for any $z\in U_{x}\backslash U_{x,1}$. Thus,
$\Upsilon_{\eta_{x},z}e=\Upsilon(\frac{e}{\eta_{x}-ze})=0$ for all
$z\in U_{x,1}$ and $\Upsilon(\frac{e}{\eta_{x}-ze})\ne 0$ for any
$z\in U_{x,2}$. This yields {\bf vi.}2.

Finally, suppose that $w\in C^{\infty}(M;\mathbb{C})$ is
holomorphic in $M$ and $\zeta:=w|_{\Gamma}$. If $z\in U_{x,1}$ and
$c\in\mathbb{C}$, then the function $\frac{w-c}{w_{x}-z}$ is
holomorphic in $M$. Hence, any $c\in\mathbb{C}$ is a solution of
(\ref{boundeq}). Now, suppose that $z\in U_{x,2}$. Since
$\frac{1}{w_{x}-z}$ has a simple pole at the point $w_{x}^{-1}(z)$
and no other poles on $M$, the function $\frac{w-c}{w_{x}-z}$ is
holomorphic in $M$ if and only if $c=w(w_{x}^{-1}(z))$. So,
(\ref{boundeq}) has a unique solution $c=w(w_{x}^{-1}(z))$ for any
$z\in U_{x,1}$. This proves {\bf vi.}3.
\smallskip

\noindent$\bullet$\,\,\,
Suppose that $\zeta,1/\zeta\in {\rm Ker}\Upsilon$. Then
$\zeta=w|_{\Gamma}$, where $w,1/w$ are holomorphic functions in
$M$. Let $U$ be an arbitrary simply connected neighborhood in $M$.
Since $w$ have no zeroes in $U$, each branch of ${\rm log} w$ is
holomorphic function in $U$. In particular, ${\rm log}
|w|={\Re\,}{\rm log\,} w$ is harmonic in $U$. Also, ${\rm log}|w|$
is single-valued function on the whole $M$. Then ${\rm
log}|w|=u^{{\rm log}|\zeta|}$ is a solution of (\ref{Eq 1}),
(\ref{Eq 2}) with $f={\rm log\,}|\zeta|$. Hence,
$\partial_{\nu}{\rm log\,}|w|=\Lambda {\rm log\,}|\zeta|$ on
$\Gamma$. Now choose $U$ in such a way that
$\overline{U}\cap\Gamma$ is a segment $\Gamma'\subset\Gamma$ of
nonzero length. Since each branch of ${\rm log\,}w$ is holomorphic
in $U$ and smooth up to $\Gamma'$, from Cauchy-Riemann conditions
it follows that
 $$
\partial_{\nu}{\rm log\,}|w|=\partial_{\nu}\Re\,{\rm log\,}w=\partial_{\gamma}\,\Im\,{\rm log\,}w=\partial_{\gamma}{\rm arg\,}w=\partial_{\gamma}{\rm arg\,}\zeta
 $$
on $\Gamma$. Therefore, $\Lambda {\rm
log\,}|\zeta|=\partial_{\gamma}{\rm arg\,}\zeta$, i.e., {\bf vii}
does hold.
\smallskip

{\it The necessity is proved}.

\subsubsection*{3. Sufficiency}

Here we assume that $\Lambda$ obeys i.-vii. and construct a
Riemannian surface $(M,g)$ such that its DN-map is $\Lambda$,
i.e., $\Lambda=\Lambda_g$ holds. Before that, we recall the known
facts and definitions that will be used in the construction.
\smallskip

\noindent$\bullet$\,\,\,A {\it commutative Banach algebra} is a
(complex) Banach space $(\mathfrak{A},\|\cdot\|)$ equipped with
the multiplication operation $\mathfrak{A}\times \mathfrak{A}\ni
\eta,\zeta\mapsto \eta\zeta\in \mathfrak{A}$ satisfying
$\eta\zeta=\zeta\eta$, $\|\eta\zeta\|\le \|\eta\|\|\zeta\|$ for
all $\eta,\zeta\in \mathfrak{A}$. Algebra $\mathfrak{A}$ is {\it
unital} if there exists $e\in\mathfrak{A}$ such that $e\eta=\eta$
holds for all $\eta\in\mathfrak{A}$. Element $\eta\in\mathfrak{A}$
is {\it invertible} if there exists $\eta^{-1}\in\mathfrak{A}$
such that $\eta^{-1}\eta=e$. The set of all $z\in\mathbb{C}$ for
which $\eta-ze$ is noninvertible is called the {\it spectrum} of
$\eta$ and is denoted by ${\rm Sp}_{\mathfrak{A}}\eta$, such a set
being compact.

A {\it character} of the commutative Banach algebra $\mathfrak{A}$
is a nonzero homomorphism $\chi: \mathfrak{A}\mapsto \mathbb{C}$.
Each character $\chi$ is a continuous map: one has
\begin{equation}\label{charcont}
|\chi(\eta)|\leqslant \|\eta\|, \qquad \eta\in\mathfrak{A}.
\end{equation}
The set of characters $\widehat{\mathfrak{A}}$ is called the {\it
spectrum} of algebra $\mathfrak{A}$. For an $\eta\in\mathfrak{A}$,
its {\it Gelfand transform} $\hat{\eta}: \
\widehat{\mathfrak{A}}\mapsto \mathbb{C}$ is defined as
$$\hat{\eta}(\chi):=\chi(\eta),\qquad\chi\in\widehat{\mathfrak{A}}.$$
For any $\hat{\eta}$, the image
$\hat{\eta}(\widehat{\mathfrak{A}})\subset\mathbb C$ coincides
with the spectrum ${\rm Sp}_{\mathfrak{A}}\eta$.

Spectrum $\widehat{\mathfrak{A}}$ is endowed with the canonical
Gelfand ($*$-weak) topology, with respect to which it is a compact
Hausdorff space. The Gelfand transforms
$\{\hat{\eta}\,|\,\,\eta\in\mathfrak{A}\}$ constitute a subalgebra
in $C(\widehat{\mathfrak{A}})$, which separates points of
$\widehat{\mathfrak{A}}$. The space $\widehat{\mathfrak{A}}$ is
connected if and only if there is no nontrivial idempotents
$\eta=\eta^{2}$, $\eta\ne 0,e$ in $\mathfrak{A}$.

A closed subset $B\subset \widehat{\mathfrak{A}}$ is called a {\it
boundary} of $\mathfrak{A}$ if
$\max_B|\hat{\eta}|=\max_{\widehat{\mathfrak{A}}}|\hat{\eta}|$ for
any $\eta\in \mathfrak{A}$. The intersection of all boundaries is
called {\it the Shilov boundary} of $\mathfrak{A}$ and denoted by
$\mathfrak{bA}$.

The key fact that we use in the proof of sufficiency is the
fundamental Bishop-Aupetit-Wermer analytic structure theorem: see
Theorem 2.2, \cite{Aupetit} or Chapter 11, \cite{Wermer}. 
\begin{Theorem}\label{Bishop theorem}
Let $\eta\in \mathfrak{A}$, the set
$\hat{\eta}(\widehat{\mathfrak{A}})\backslash\hat{\eta}(\mathfrak{bA})$
is non-empty, and $V$ is its connected component. Suppose that the
set $\{z\in V \ | \ \sharp\,\hat{\eta}^{-1}(\{z\})<\infty\}$ is of
nonzero Lebesgue measure. Then $\sharp\,\hat{\eta}^{-1}(\{z\})\le
N<\infty$ for any $z\in V$ and the subset
$\hat{\eta}^{-1}(V)\subset \widehat{\mathfrak{A}}$ has the
structure of 1-dim complex analytic manifold, on which all
functions $\hat{\zeta}$ ($\zeta\in\mathfrak{A}$) are holomorphic.
\end{Theorem}
The rest of the proof of Theorem \ref{main} is as follows. We
construct a Riemann surface $M$ as the spectrum
$\widehat{\mathfrak{A}}$ of some Banach function algebra
$\mathfrak{A}$ provided by conditions {\bf i} and {\bf iii}. Then,
using Theorem \ref{Bishop theorem} and the condition {\bf v}, we
endow a part $\Omega_{\eta}\subset\widehat{\mathfrak{A}}$ with the
structure of Riemannian surface, and this part depends on the
element $\eta\in \mathfrak{A}$. The condition {\bf vi} enables
one, by varying the elements $\eta\in\mathfrak{A}$, to cover the
whole spectrum $\widehat{\mathfrak{A}}$ by its analytic parts
$\Omega_{\eta}$ and thus endow the $\widehat{\mathfrak{A}}$ with
the structure of Riemannian surface. Also, due to Theorem
\ref{Bishop theorem}, the Gelfand transforms of elements of
$\mathfrak{A}$ form a subalgebra in the algebra of holomorphic
smooth functions on $M$. By conditions {\bf ii} and {\bf iv}, this
subalgebra coincides with the set of all holomorphic smooth
functions on $M$. By the latter, we show that $\Lambda$ coincides
with the DN-map of the surface $M$ on all traces on $\Gamma$ of
real parts of holomorphic smooth functions on $M$, the set of such
traces being of finite codimension. To check that $\Lambda$
coincides with the DN-map of $M$ on {\it all other functions} from
$C^{\infty}(\Gamma;\mathbb{R})$, we use the remaining condition
{\bf vii}.

So, we proceed to prove the sufficiency of the conditions of
Theorem {\rm \ref{main}}.
\smallskip

\noindent$\bullet$\,\,\,In view of {\bf i}, the set ${\rm
Ker}\Upsilon$ is a unital (sub)algebra in $C(\Gamma;\mathbb{C})$.
The closure $$\mathfrak{A}\,:=\,\overline{{\rm
Ker}\Upsilon}\,\subset C(\Gamma;\mathbb{C})$$ is a unital
commutative Banach algebra with the norm
$\parallel\zeta\parallel:=\max_{\Gamma}|\zeta|$. We denote its
spectrum $\widehat{\mathfrak{A}}$ by $M$. Recall that $M$ is a
Hausdorff compact space. Also, $M$ is connected: indeed, if
$\eta^{2}=\eta$ on $\Gamma$, then $\eta(x)=0$ or $1$ for any
$x\in\Gamma$; since $\eta$ is continuous and $\Gamma$ is
connected, this means that $\eta=0$ or $1$. Note that the set of
smooth elements of $\mathfrak{A}$ coincides with ${\rm
Ker}\Upsilon$ due to property {\bf iii}.

The Dirac measures $\delta_{x} : \ \mathfrak{A}\ni \zeta\mapsto
\zeta(x)\in\mathbb{C}$, $x\in\Gamma$ constitute a subset
$\delta_{\Gamma}\subset M$. In view of (\ref{charcont}) and the
definition of the norm $\|\cdot\|$, one has
$|\hat{\eta}(\chi)|=|\chi(\eta)|\leqslant\|\eta\|=\max_{x\in\Gamma}|\delta_{x}(\eta)|$
for any $\chi\in M$ and $\eta\in\mathfrak{A}$. Hence,
$\delta_{\Gamma}$ is a boundary of $\mathfrak{A}$ and thus it
contains the Shilov boundary $\mathfrak{bA}$ of $\mathfrak{A}$.
\smallskip

\noindent$\bullet$\,\,\,Our first goal is to endow
$M\backslash\delta_{\Gamma}$ with the structure of an analytic
manifold via Theorem \ref{Bishop theorem}. To verify the
conditions of Theorem \ref{Bishop theorem}, we prove that
$\hat{\eta}^{-1}(\{z\})$ is finite for any $\eta\in {\rm
Ker}\Upsilon=\mathfrak{A}\cap C^{\infty}(\Gamma;\mathbb{C})$ and
$z\in \hat{\eta}(\widehat{\mathfrak{A}})\backslash\eta(\Gamma)$.
The proof is based on a bijection between the characters from
$\hat{\eta}^{-1}(\{z\})$ and the characters over a certain
finite-dimensional factor-algebra $\mathfrak{A}_{\eta,z}$ which is
constructed below. Let's get down to implementing this plan.
\medskip

Let $\eta\in {\rm Ker}\Upsilon$ and
$z\in\mathbb{C}\backslash\eta(\Gamma)$; then the function
$\eta-ze$ is invertible in $C^{\infty}(\Gamma;\mathbb{C})$ (but
not necessarily in $\mathfrak{A}$). Consider the main ideal
$\mathcal{I}_{\eta,z}:=(\eta-ze)\,\mathfrak{A}$ in $\mathfrak{A}$.
It is closed in $\mathfrak{A}$: indeed, if
$\mathcal{I}_{\eta,z}\ni\zeta_{k}\to\zeta$ in
$C(\Gamma;\mathbb{C})$, then the convergence
$\mathfrak{A}\ni\frac{\zeta_{k}}{\eta-ze}\to\frac{\zeta}{\eta-ze}$
holds by $\frac{1}{\eta-ze}\in C^{\infty}(\Gamma;\mathbb{C})$.
Since $\mathfrak{A}$ is Banach, $\frac{\zeta}{\eta-ze}\in
\mathfrak{A}$ and $\zeta\in \mathcal{I}_{\eta,z}$. Since ${\rm
Ker}\Upsilon$ is dense in $\mathfrak{A}$, the set
$\mathcal{I}_{\eta,z}^{\infty}:=(\eta-ze)\,{\rm Ker}\Upsilon$ is
dense in $\mathcal{I}_{\eta,z}$. The function $\zeta\in {\rm
Ker}\Upsilon$ belongs to $\mathcal{I}_{\eta,z}^{\infty}$ if and
only if
$0=\Upsilon(\frac{\zeta}{\eta-ze})=\Upsilon_{\eta,z}(\zeta)$.

Introduce the factor-algebra
$$\mathfrak{A}_{\eta,z}:=\mathfrak{A}/\mathcal{I}_{\eta,z}$$ with
the factor-norm
$\|\zeta+\mathcal{I}_{\eta,z}\|_{\eta,z}:=\inf_{\tilde{\zeta}\in\mathcal{I}_{\eta,z}}\|\zeta+\tilde{\zeta}\|$;
here and in what follows we denote by $\zeta+\mathcal{I}_{\eta,z}$
the equivalence class in $\mathfrak{A}_{\eta,z}$ of element
$\zeta\in\mathfrak{A}$. Due to definition of the factor-norm and
the equality $\overline{{\rm Ker}\Upsilon}=\mathfrak{A}$, the set
$\mathfrak{A}_{\eta,z}^{\infty}:=\{\zeta+\mathcal{I}_{\eta,z} \ |
\ \zeta\in {\rm Ker}\Upsilon\}$ is dense in
$\mathfrak{A}_{\eta,z}$. Let us prove that the algebra
$\mathfrak{A}_{\eta,z}$ is finite-dimensional. To this end,
consider a linear map $\mathscr{G}_{\eta,z}: \
\mathfrak{A}_{\eta,z}^{\infty}\mapsto
C^{\infty}(\Gamma;\mathbb{C})$ defined by the rule
$$\mathscr{G}_{\eta,z}(\zeta+\mathcal{I}_{\eta,z})=\Upsilon_{\eta,z}(\zeta).$$
The map $\mathscr{G}_{\eta,z}$ is well-defined and its kernel is
trivial. Indeed, if
$\zeta_{1}+\mathcal{I}_{\eta,z}=\zeta_{2}+\mathcal{I}_{\eta,z}\in
\mathfrak{A}_{\eta,z}^{\infty}$, then $\zeta_{1}-\zeta_{2}\in
\mathcal{I}_{\eta,z}\cap
C^{\infty}(\Gamma;\mathbb{C})=\mathcal{I}_{\eta,z}^{\infty}$ and
$\Upsilon_{\eta,z}(\zeta_{1})-\Upsilon_{\eta,z}(\zeta_{2})=0$.
Similarly, if
$\mathscr{G}_{\eta,z}(\zeta+\mathcal{I}_{\eta,z})=0$, then
$\Upsilon_{\eta,z}(\zeta)=0$ and thus $\zeta\in
\mathcal{I}_{\eta,z}^{\infty}\subset\mathcal{I}_{\eta,z}$ i.e.
$\zeta+\mathcal{I}_{\eta,z}$ is the zero element in
$\mathfrak{A}_{\eta,z}$. Note that
$\mathscr{G}_{\eta,z}\mathfrak{A}_{\eta,z}^{\infty}=\Upsilon_{\eta,z}{\rm
Ker}\Upsilon$. Since the map $\mathscr{G}_{\eta,z}$ is a bijection
of $\mathfrak{A}_{\eta,z}^{\infty}$ and
$\mathscr{G}_{\eta,z}\mathfrak{A}_{\eta,z}^{\infty}$, one has
$${\rm dim\,}\mathfrak{A}_{\eta,z}^{\infty}={\rm dim\,}[\Upsilon_{\eta,z}{\rm Ker}\Upsilon].$$
In view of condition {\bf v}, the right-hand side is equal to the
integral
$$\frac{1}{2\pi i}\int\limits_{\Gamma}\frac{\partial_{\gamma}\eta}{\eta-ze}\,d\gamma\,.$$
Since the functions $\eta$ and $\frac{1}{\eta-ze}$ are smooth,
this integral is finite. So, ${\rm
dim\,}\mathfrak{A}_{\eta,z}^{\infty}$ is finite and, since
$\mathfrak{A}_{\eta,z}^{\infty}$ is dense in
$\mathfrak{A}_{\eta,z}$, one has
$${\rm dim\,}\mathfrak{A}_{\eta,z}=\frac{1}{2\pi i}\int\limits_{\Gamma}\frac{\partial_{\gamma}\eta}{\eta-ze}\,d\gamma\,.$$
Note that the right-hand side is the number $d(z)$ of revolutions
of the image $\eta(\Gamma)\subset\mathbb C$  
around the point $z$; this number depends only on the connected
component $V$ of $\mathbb{C}\backslash\eta(\Gamma)$ that contains
$z$. If $d(z)=0$, then ${\rm dim\,}\mathfrak{A}_{\eta,z}=0$ and
$e\in\mathcal{I}_{\eta,z}=\mathfrak{A}$. This means that $\eta-ze$
is invertible in $\mathfrak{A}$ and $z\not\in{\rm
Sp}_{\mathfrak{A}}\eta=\hat{\eta}(M)$. Thus,
$$\hat{\eta}(M)\backslash\eta(\Gamma)=\{z\in\mathbb{C}\backslash\eta(\Gamma) \ | \ d(z)>0\}.$$
\smallskip

\noindent$\bullet$\,\,\,
Now, we show that the set $\hat{\eta}^{-1}(\{z\})$ is finite,
$z\in V$ being the same as before. Let $\tilde{\chi}$ be a
character on the algebra $\mathfrak{A}_{\eta,z}$; then the rule
\begin{equation}
\label{factrochar}
\chi(\zeta):=:\tilde{\chi}(\zeta+\mathcal{I}_{\eta,z})
\end{equation}
defines a character $\chi\in M$ that vanishes on
$\mathcal{I}_{\eta,z}$. Hence, we have $\chi(\eta-ze)=0$ and
$\chi(\eta)=z$. Conversely, suppose that $\chi\in M$ and
$\chi(\eta)=z$ (then, obviously,
$\chi(\mathcal{I}_{\eta,z})=\{0\}$). Then the same rule
(\ref{factrochar}) defines the character $\tilde{\chi}$ on
$\mathfrak{A}_{\eta,z}$. Thus, we have
\begin{equation}
\label{sharp} {\sharp\,}
\widehat{\mathfrak{A}}_{\eta,z}={\sharp\,}\hat{\eta}^{-1}(\{z\}).
\end{equation}
Suppose that $\tilde{\chi}_{1},\dots,\tilde{\chi}_{N}$ are the
distinct characters in $\widehat{\mathfrak{A}}_{\eta,z}$. Since
the Gelfand transforms of the elements
$\zeta+\mathcal{I}_{\eta,z}\in\mathfrak{A}_{\eta,z}$ separate the
points of $\widehat{\mathfrak{A}}_{\eta,z}$, there exists
$\zeta_{ij}\in\mathfrak{A}$ such that
$\tilde{\chi}_{i}(\zeta_{ij}+\mathcal{I}_{\eta,z})=1$ and
$\tilde{\chi}_{j}(\zeta_{ij}+\mathcal{I}_{\eta,z})=0$. Denote
$q_{i}:=\Pi_{j\ne i}(\zeta_{ij}+\mathcal{I}_{\eta,z})$, then
$\tilde{\chi}_{j}(q_{i})=\delta_{ij}$. In particular, $q_{i}$,
$i=1,\dots,N$ are linearly independent in $\mathfrak{A}_{\eta,z}$.
Therefore, $N\leqslant {\rm dim\,}\mathfrak{A}_{\eta,z}=d(z)$. So,
we see that ${\sharp\,} \widehat{\mathfrak{A}}_{\eta,z}\leqslant
d(z)$ and, by (\ref{sharp}), arrive at
${\sharp\,}\hat{\eta}^{-1}(\{z\})\leqslant d(z)<\infty$.
\smallskip

\noindent$\bullet$\,\,\,
Suppose that $z\in \hat{\eta}(M)\backslash\eta(\Gamma)$ (to
provide $\hat{\eta}(M)\backslash\eta(\Gamma)\ne\varnothing$, one
can take $\eta=\eta_{x}$ for any function $\eta_{x}$ obeying {\bf
vi.}1). Denote by $V$ the connected component of
$\mathbb{C}\backslash\eta(\Gamma)$ that contains $z$. Obviously,
$V$ is open and, hence, it has a nonzero Lebesgue measure.
Moreover, $1\leqslant d(z)<\infty$ and $d(z')=d(z)$ is valid for
any $z'\in V$. For such a $z'$, the dimension $d(z')$ of algebra
$\mathfrak A_{\eta,z'}$ is finite and nonzero. Hence,
$\mathcal I_{\eta,z'}\not=\mathfrak A$, i.e., $\eta-z'e$ is
noninvertible. Therefore $z'\in\hat\eta(M)$ and, thus, for the
whole $V$ the embedding
$V\subset\hat{\eta}(M)\backslash\eta(\Gamma)$ holds.

Since $\eta(\Gamma)=\hat{\eta}(\delta_{\Gamma})$ and
{$\mathfrak{bA}\subset \delta_{\Gamma}$}, 
the set $V$ does not intersect with $\hat{\eta}(\mathfrak{bA})$.
In view of Theorem \ref{Bishop theorem}, the set
$\hat{\eta}^{-1}(V)\subset \widehat{\mathfrak{A}}$ has the
structure of 1-dim complex analytic manifold on which all
functions $\hat{\zeta}$ ($\zeta\in\mathfrak{A}$) are analytic.
Thus, for any character $\chi\in \hat{\eta}^{-1}(V)$ there
exist an open (in the Gelfand topology) neighborhood $U\ni\chi$
and a homeomorphism $\kappa: U\rightarrow D$ onto an open disk
$D\subset\mathbb{C}$ such that any function
$\hat{\zeta}\circ\kappa^{-1}$ ($\zeta\in\mathfrak{A}$) is
holomorphic on $D$. In other words, every $\chi\in
\hat{\eta}^{-1}(V)$ does possess a local analytic coordinate
$\hat\eta$.

As a result, we can represent the spectrum $M$ as the disjoint
union:
$$M=M'\cup\delta_{\Gamma}\cup\widetilde{M},$$
where
$$M':=\bigcup_{\eta\in{\rm Ker}\Upsilon}\hat{\eta}^{-1}(\mathbb{C}\backslash\eta(\Gamma))$$
is the set of characters that can be provided with the local
coordinate by the choice of a suitable $\eta\in {\rm
Ker}\Upsilon$, and $\widetilde{M}:=M\backslash
(M'\cup\delta_{\Gamma})$.
\smallskip

\noindent$\bullet$\,\,\,
Let us show that $\widetilde{M}=\varnothing$. Suppose, on the
contrary, that $\chi\in\widetilde{M}$ and $\eta$ satisfies
condition {\bf vi.}1 (as such $\eta$, one can choose any
$\eta_{x}$ from {\bf vi}). Then, $z:=\hat{\eta}(\chi)\in
\eta(\Gamma)$ and the set
$\delta_{\Gamma}\cap\hat{\eta}^{-1}(\{z\})=\{\delta_{x} \ | \
\eta(x)=z\}$ is finite. Denote by
$\delta_{x_{1}},\dots,\delta_{x_{l}}$ all characters from
$\delta_{\Gamma}\cap\hat{\eta}^{-1}(\{z\})$. Choose $\eta_{x_{k}}$
from condition {\bf vi} in such a way that $\eta_{x_{k}}(x_{k})=z$
(this condition can always be satisfied since  $\eta_{x}$ in {\bf
vi} are determined up to a constant). In view of {\bf vi.}3, any
$\zeta\in{\rm Ker}\Upsilon$ can be represented as
$\zeta=c_{\zeta,k}e+\zeta'_{k}$, where
$\Upsilon_{\eta_{x_{k}},z}(\zeta'_{k})=0$. Then
$\tilde{\zeta}'_{k}:=\frac{\zeta'_{k}}{\eta_{x_{k}}-ze}$ belongs
to ${\rm Ker}\Upsilon=\mathfrak{A}\cap
C^{\infty}(\Gamma;\mathbb{C})$, whence $\zeta'_{k}(x_k)=0$ and
$c_{\zeta,k}=\zeta(x_{k})$. Thus,
$\zeta=\zeta(x_{k})e+(\eta_{x_{k}}-ze)\tilde{\zeta}'_{k}$ with
$\tilde{\zeta}'_{k}\in {\rm Ker}\Upsilon$ and one has
 $$
\chi(\zeta)-\zeta(x_{k})=[\chi(\eta_{x_{k}})-z]\,\chi(\tilde{\zeta}'_{k}).
 $$
In particular, if $\chi(\eta_{x_{k}})=z$, then $\chi$ coincides
with $\delta_{x_{k}}$ on the dense set ${\rm Ker}\Upsilon$ in
$\mathfrak{A}$ and, hence, on the whole $\mathfrak{A}$. Thus, the
embedding $\chi\in\widetilde{M}$ shows that $\chi(\eta_{x_{k}})\ne
z$ for any $k=1,\dots,l$. Then the element
$$\theta=\prod\limits_{k=1}^{l}[\,\eta_{x_{k}}-\eta_{x_{k}}(x_{k})\,e\,]$$
satisfies $\chi(\theta)\ne 0$ and $\theta(x_{k})=0$ for any
$k=1,\dots,l$.

Denote
$$\eta_{\varepsilon,\varphi}:=\eta+\varepsilon^{2}e^{i\varphi}\chi(\theta)^{-1}\theta,$$
where $\varepsilon>0$ and $\varphi\in [0,2\pi)$. Due to condition
{\bf vi.}1, the function $\eta-ze$ has a zero of multiplicity
$\leqslant m<\infty$ at each $x_{k}$ and there are no other zeros
of $\eta-ze$ on $\Gamma$. Thus, the pre-image
$\eta^{-1}(B_{\varepsilon})$ of the $\varepsilon-$neighborhood
$U_{\varepsilon}$ of $z$ is contained in
$O(\varepsilon^{1/m})$-neighborhood of the set
$\{x_{1},\dots,x_{M}\}$ in $\Gamma$. Therefore,
\begin{align}
\notag &
|\eta_{\varepsilon,\varphi}(x)-\eta(x)|=\varepsilon^{2}|\chi(\theta)^{-1}\theta(x)|=\varepsilon^{2}|\theta(x_{k})+O({\rm
dist}\{x,x_{k}\})|=\\
\label{est1} &
=\varepsilon^{2}|0+O(\varepsilon^{1/m})|=O(\varepsilon^{2+1/m}),
\end{align}
where $x_{k}$ is chosen to be the nearest point to $x$. For $x\in
\Gamma\backslash\eta^{-1}(B_{\varepsilon})$, one has
$|\eta(x)|\geqslant\varepsilon$, whence
\begin{equation}
\label{est2}
|\eta_{\varepsilon,\varphi}(x)|=|\eta(x)+\varepsilon^{2}O(1)|\geqslant
\varepsilon-O(\varepsilon^{2}).
\end{equation}
Estimates (\ref{est1}),(\ref{est2}) show that, for sufficiently
small $\varepsilon$, the set $B_{\varepsilon/2}$ does not
intersect with the fragment
$\eta_{\varepsilon,\varphi}(\Gamma\backslash\eta^{-1}(B_{\varepsilon}))$
of $\eta_{\varepsilon,\varphi}(\Gamma)$ while the fragment
$\eta_{\varepsilon,\varphi}(\eta^{-1}(B_{\varepsilon}))$ is
contained in $O(\varepsilon^{2+1/m})-$neighbourhood of
$\eta(\eta^{-1}(B_{\varepsilon}))$. Thus, it is possible to choose
$\varepsilon$ and $\varphi$ such that
$\hat{\eta}_{\varepsilon,\varphi}(\chi)=\chi(\eta_{\varepsilon,\varphi})=z+\varepsilon^{2}e^{i\varphi}\not\in
\eta_{\varepsilon,\varphi}(\Gamma)$. This means that
$\chi\in\hat{\eta}_{\varepsilon,\varphi}^{-1}(\mathbb{C}\backslash\eta(\Gamma))\subset
M'$, so that we arrive at the contradiction and prove that
$\widetilde{M}=\varnothing$.
\smallskip

\noindent${\bullet}$\,\,\,
Now, we coordinatize $\delta_{\Gamma}$. Let $\delta_{x}$
($x\in\Gamma$) be an arbitrary character from $\delta_{\Gamma}$.
Consider the map
$$\hat{\eta}_{x}: \ \hat{\eta}_{x}^{-1}(U_{x})\mapsto\mathbb{C},$$
where $\eta_{x}$, $U_{x}$ are the same as in condition {vi.} Due
to condition {\bf vi.}2, for any $z$ from one connected component
$U_{x,1}$ of $U_{x}\backslash\eta_{x}(\Gamma)$ one has
$\frac{1}{\eta_{x}-ze}\in{\rm Ker}\Upsilon\subset\mathfrak{A}$.
Hence $z\not\in{\rm Sp}_{\mathfrak{A}}\eta_{x}=\hat{\eta}_{x}(M)$.
Now, suppose that $z\in U_{x}\backslash U_{x,1}$ and $\zeta\in{\rm
Ker}\Upsilon$. By {\bf vi.}2, one has
$\frac{1}{\eta_{x}-ze}\not\in{\rm Ker}\Upsilon$, i.e., either
$z\in\eta_{x}(\Gamma)$ or $\frac{1}{\eta_{x}-ze}$ is smooth on
$\Gamma$ and does not belong to $\mathfrak{A}$. Thus, $z\in{\rm
Sp}_{\mathfrak{A}}\eta_{x}=\hat{\eta}_{x}(M)$. Also, due to {\bf
vi.}3, there exists $c_{\zeta,z}\in\mathbb{C}$ such that
$\frac{\zeta-c_{\zeta,z}e}{\eta_{x}-ze}=:\tilde{\zeta}_{\zeta,z}\in
{\rm Ker}\Upsilon\subset\mathfrak{A}$. If
$\chi\in\hat{\eta}_{x}^{-1}(\{z\})$, then one has
 \begin{align*}
&\chi(\zeta)=\chi(c_{\zeta,z}e+(\eta_{x}-ze)\tilde{\zeta}_{\zeta,z})=c_{\zeta,z}+(\chi(\eta_{x})-z)\chi(\tilde{\zeta}_{\zeta,z})=\\
&=c_{\zeta,z}+(\hat{\eta}_{x}(\chi)-z)\chi(\tilde{\zeta}_{\zeta,z})=c_{\zeta,z}.
 \end{align*}
Since $\zeta\in{\rm Ker}\Upsilon$ is arbitrary and ${\rm
Ker}\Upsilon$ is dense in $\mathfrak{A}$, this means that
$\chi(\zeta)=\chi'(\zeta)$ for any
$\chi,\chi'\in\hat{\eta}_{x}^{-1}(z)$ and $\zeta\in\mathfrak{A}$.
Thus, ${\sharp\,}\hat{\eta}_{x}^{-1}(\{z\})=1$ and the map
$\hat{\eta}_{x}: \ \hat{\eta}_{x}^{-1}(U_{x})\mapsto
U_{x}\backslash U_{x,1}$ is a bijection.

Since the expansion
 $$
\zeta=c_{\zeta,z}e+(\eta_{x}-ze)\tilde{\zeta}_{\zeta,z}
\,\,\,\text{with}\,\,\,
c_{\zeta,z}=\hat{\zeta}\circ\hat{\eta}_{x}^{-1}(z)=
\hat{\zeta}\left(\hat{\eta}_{x}^{-1}(z)\right)\,\,\,\text{and}\,\,\,
\tilde{\zeta}_{\zeta,z}\in{\rm Ker}\Upsilon
 $$
holds for any $z\in U_{x}\backslash U_{x,1}$ and $\zeta\in{\rm
Ker}\Upsilon$, one can iterate it to obtain
\begin{align*}
&
\zeta=c_{\zeta,z}e+(\eta_{x}-ze)\left[\,c'_{\zeta,z}e+(\eta_{x}-ze)\tilde{\zeta}'_{\zeta,z}\right]\,\,\,\text{with}\,\,\,
c_{\zeta,z}=\hat{\zeta}\circ\hat{\eta}_{x}^{-1}(z), \\
&
c'_{\zeta,z}=\hat{\tilde{\zeta}}_{\zeta,z}\circ\hat{\eta}_{x}^{-1}(z)\,\,\,\text{and}\,\,\,
\tilde{\zeta}'_{\zeta,z}\in{\rm Ker}\Upsilon
 \end{align*}
for any $z\in U_{x}\backslash U_{x,1}$. Then one has
\begin{align}
\notag &
\hat{\zeta}\circ\hat{\eta}_{x}^{-1}(z')-\hat{\zeta}\circ\hat{\eta}_{x}^{-1}(z)=\\
\notag & =
[\hat{\eta}_{x}^{-1}(z')]\left[\hat{\zeta}\circ\hat{\eta}_{x}^{-1}(z)e+(\eta_{x}-ze)
(\hat{\tilde{\zeta}}_{\zeta,z}\circ\hat{\eta}_{x}^{-1}(z)e+(\eta_{x}-ze)\tilde{\zeta}'_{\zeta,z})\right]-\\
\label{holocontinuity} &
-\hat{\zeta}\circ\hat{\eta}_{x}^{-1}(z)=(z'-z)\left[\hat{\tilde{\zeta}}_{\zeta,z}\circ\hat{\eta}_{x}^{-1}(z)+(z'-z)\cdot\hat{\tilde{\zeta}}'_{\zeta,z}\circ\hat{\eta}_{x}^{-1}(z')\right]
\end{align}
for any $z,z'\in U_{x}\backslash U_{x,1}$. In view of
(\ref{charcont}), $$ |{\hat
\zeta'}\circ\hat{\eta}_{x}^{-1}(z)|=|[\hat{\eta}_{x}^{-1}(z)]({\zeta'})|\leqslant\|{
\zeta'}\|, \quad \zeta'\in\mathfrak{A}.$$ Then the function
$z'\mapsto\hat{\tilde{\zeta}}'_{\zeta,z}\circ\hat{\eta}_{x}^{-1}(z')$
is bounded on $U_{x}\backslash U_{x,1}$ and (\ref{holocontinuity})
implies that the function
$z\mapsto\hat{\zeta}\circ\hat{\eta}_{x}^{-1}(z)$ is continuous on
$U_{x}\backslash U_{x,1}$ for any $\zeta\in{\rm Ker}\Upsilon$. The
same is true for the function
$z'\mapsto\hat{\tilde{\zeta}}'_{\zeta,z}\circ\hat{\eta}_{x}^{-1}(z')$.
Then (\ref{holocontinuity}) implies also that, for any
$\zeta\in{\rm Ker}\Upsilon$, the function
$z\mapsto\hat{\zeta}\circ\hat{\eta}_{x}^{-1}(z)$ is holomorphic on
$U_{x}\backslash\overline{U_{x,1}}$ and its partial derivatives
are continuous on $U_{x}\backslash U_{x,1}$. In view of definition
of the Gelfand topology, the map $\eta_{x}: \
\hat{\eta}_{x}^{-1}(U_{x})\mapsto U_{x}\backslash U_{x,1}$ is
homeomorphism. So, any character $\delta_{x}\in\delta_{\Gamma}$ is
coortinatizable in the following sense: there exists a
neighbourhood $V:=\hat{\eta}_{x}^{-1}(U_{x})$ (in the Gelfand
topology) of $\delta_{x}$ and the local coordinate
$\hat{\eta}_{x}: \ V\mapsto U_{x}\backslash U_{x,1}$ in which all
functions $\hat{\zeta}$ ($\zeta\in{\rm Ker}\Upsilon$) are
holomorphic on $U_{x}\backslash\overline{U_{x,1}}$ and continuous
differentiable up to the (smooth) curve
$U_{x}\cap\eta_{x}(\Gamma)$. Note that, in view of vi, a., the map
$\eta^{-1}(U_{x})\ni
x'\rightarrow(\Re\eta_{x}(\delta_{x'}),\Im\eta_{x}(\delta_{x'})$
is a diffeomorphism.
\smallskip

\noindent${\bullet}$\,\,\,
We have proved above that the whole
$M$ is coordinatizable i.e. each character $\chi\in M$ has the
open (in the Gelfand topology) neighbourhood $V_{\chi}$ and the
homeomorphism $\kappa_{\chi}:\ V_{\chi}\mapsto
U_{\chi}\subset\mathbb{C}$ such that 1) the set
$U'_{\chi}:=\kappa_{\chi}(V_{\chi}\backslash\chi(\delta_{\Gamma}))$
is open, 2) $U_{\chi}\backslash U'_{\chi}\subset\partial U_{\chi}$
is empty or it is the fragment of smooth curve, and 3) each
function $\hat{\zeta}\circ\kappa_{\chi}^{-1}$ ($\zeta\in{\rm
Ker}\Upsilon$) is holomorphic on $U'_{\chi}$ and continuous
differentiable up to $U'_{\chi}\subset\partial U_{\chi}$.

Now, we construct a biholomorphic atlas on $M$ using
$\{V_{\chi},\kappa_{\chi}\}_{\chi\in M}$. The collection
$\{V_{\chi}\}_{\chi\in M}$ is an open cover of $M$ and, since $M$
is compact, there exists a finite subcover
$\{V_{\chi_{k}}\}_{k=1}^{L}$. Denote $V_{k}:=V_{\chi_{k}}$ and
$\kappa_{k}:=\kappa_{\chi_{k}}$. Suppose that $V_{k}\cap
V_{l}\ne\varnothing$ and denote $W_{k}:=\kappa_{k}((V_{k}\cap
V_{l})\backslash\delta_{\Gamma})$, $W_{l}:=\kappa_{l}((V_{k}\cap
V_{l})\backslash\delta_{\Gamma})$. Choose an arbitrary nonconstant
$\zeta\in{\rm Ker}\Upsilon$ (for example, one of $\eta_{x}$ from
condition iii.), then $\hat{\zeta}\circ\kappa_{k}^{-1}$,
$\hat{\zeta}\circ\kappa_{l}^{-1}$ are holomorphic on $W_{k}$ and
$W_{l}$, respectively. In particular, any zero of $\nabla
\Re(\hat{\zeta}\circ\kappa_{k}^{-1})$ on $W_{k}$ is isolated. If
$\kappa_{k}(\chi)\in W_{k}$ does not coincide with zero of
$\nabla\Re(\hat{\zeta}\circ\kappa_{k}^{-1})$, then there exists
the neighbourhood $V'$ of $\chi$ such that
$\zeta\circ\chi_{k}^{-1}: \ \kappa_{k}(V')\mapsto
\zeta\circ\chi_{k}^{-1}(W')$ is biholomorphic map. So, the
function
$$\kappa_{k}\circ\kappa_{l}^{-1}=\kappa_{k}\circ\zeta^{-1}\circ\zeta\circ\kappa_{l}=(\zeta\circ\kappa_{k}^{-1})^{-1}\circ(\zeta\circ\kappa_{l}^{-1})$$
is holomorphic on $\kappa_{l}(V')$. So,
$\kappa_{k}\circ\kappa_{l}^{-1}$ is holomorphic on $W_{l}$ except
for only some isolated points. Since
$\kappa_{k}\circ\kappa_{l}^{-1}$ is continuous on $W_{l}$, we find
that $\kappa_{k}\circ\kappa_{l}^{-1}$ is holomorphic on the whole
$W_{l}$. The same reason shows that
$\kappa_{l}\circ\kappa_{k}^{-1}$ is holomorphic on $W_{k}$ and,
thus, the transition function $\kappa_{k}\circ\kappa_{l}^{-1}$ is
biholomorphic. So, we have proved that
$\{V_{k}:=V_{\chi_{k}},\kappa_{k}:=\kappa_{\chi_{k}}\}_{k=1}^{L}$
is a biholomorphic atlas on $M$. Endowed with this atlas, $M$ is a
Riemann surface with boundary $\delta_{\Gamma}$. Moreover, the map
$\delta : \ \Gamma\ni x\mapsto\delta_{x}\in\delta_{\Gamma}$ is
diffeomorfism. In what follows, we identify $\Gamma$ and
$\delta_{\Gamma}$ via the map $\delta$.
\smallskip

\noindent${\bullet}$\,\,\,
Now, we introduce the metric $g$ and the rotation $\Phi$ on $M$
which are consistent with the metrics and the tangent field
$\gamma$ on $\Gamma$. Endow $M$ with the metric tensor
$g'=\sum_{k=1}^{L}\psi_{k}g_{k}$, where $g_{k}^{ij}=\delta^{ij}$
in local coordinates $\kappa_{k}$, and $\{\psi_{k}\}_{k=1}^{L}$ is
a partition of unity on $M$: $\psi_{k}\circ\kappa^{-1}_{l}$ is
smooth for any $l$, $\psi_{k}\geqslant 0$, ${\rm
supp}\psi_{k}\subset V_{k}$, and $\sum_{k=1}^{L}\psi_{k}=1$. Since
the transition functions are biholomorphic, the tensor $g'$ is of
the form
$\sum_{k=1}^{L}\psi_{k}|\nabla\Re(\kappa_{k}\circ\kappa^{-1}_{l})|^{2}\delta^{ij}$
in any local coordinates $\kappa_{k}$. Tensor $g'$ induces the new
metrics $d\gamma'=q(x)d\gamma$ on $\Gamma\equiv\delta_{\Gamma}$,
where $q>0$ is smooth on $\Gamma$ due to condition {\bf vi.}1.
Introducing a smooth conformal multiplier $\rho$, such that
$\rho=q^{-1}$ on $\Gamma$, we obtain the new metric tensor $g=\rho
g'$ which is consistent with the original metric on $\Gamma$.

Choose a continuous family of rotations $M\ni x\mapsto\Phi_{x}\in
{\rm End\,}T_{x}M$,
$$g(\Phi_{x}a,\Phi_{x}b)=g(a,b), \quad g(\Phi_{x} a,a)=0, \qquad \forall a,b\in T_{x}M, \ x\in M$$
such that $\Phi^{1}_{1}=\Phi_{2}^{2}=0$,
$\Phi_{2}^{1}=-\Phi_{1}^{2}=1$ in local coordinates
$x_{1}=\Re\kappa_{k}$, $x_{2}=\Im\kappa_{k}$. For any $k$ and
$\zeta\in{\rm Ker}\Upsilon$, the function
$\hat{\zeta}\circ\kappa_{k}$ is holomorphic, whence
\begin{equation}
\label{CR--} \nabla\Im\hat{\zeta}=\Phi\nabla\Re\hat{\zeta} \text{
in } M\backslash\Gamma.
\end{equation}
So, any function $\hat{\zeta}$ ($\zeta\in{\rm Ker}\Upsilon$) is
holomorphic in $M\backslash\Gamma$ (in the sense of Cauchi-Riemann
conditions (\ref{CR--})) and continuously differentiable up to
$\Gamma$. In particular, any functions $\Re\hat{\zeta}$ and
$\Im\hat{\zeta}$ are harmonic in $M\backslash\Gamma$ and
continuously differentiable up to $\Gamma$. Let $\nu$ be the
outward normal on $\Gamma$. Then $\gamma'=\Phi\nu$ is the tangent
field on $\Gamma$ and, hence, it coincides with $s\gamma$, where
$s=1$ or $-1$. Choose some $\eta\in{\rm Ker}\Upsilon$ and
$z\in{\rm Sp}_{\mathfrak{A}}\eta\backslash\eta(\Gamma)$, then
$\hat{\eta}-ze$ have at least one zero on $M$. Note that
$\hat{\eta}=\eta$ on $\Gamma$. Consider the integral
$$\int\limits_{\Gamma}\frac{1}{2\pi i}\frac{\partial_{\gamma'}\hat{\eta}}{\hat{\eta}-z}d\gamma=s\int\limits_{\Gamma}\frac{1}{2\pi i}\frac{\partial_{\gamma}\eta}{\eta-ze}d\gamma.$$
In view of the argument principle, the integral in the left-hand
side coincides with the number of zeroes of $\hat{\eta}$ counted
with their multiplicities, and, thus, it is positive. The integral
in the right-hand side is positive in view of (\ref{algebraic
argument principle}). Therefore, $s=1$ and $\gamma=\gamma'$.
\smallskip

\noindent${\bullet}$\,\,\,
Suppose that $f\in {\rm Ker}(\partial_{\gamma}+\Lambda J\Lambda)$.
Denote $h:=J\Lambda f(=J\Lambda J\partial_{\gamma}f)$ and
$\zeta=f+ih$. Then $\partial_{\gamma}h=\Lambda f$,
$\partial_{\gamma}f=-\Lambda h$ and, hence, $\zeta\in{\rm
Ker}\Upsilon$. In view of what already proven, the functions
$u:=\Re\hat{\zeta}$, $v:=\Im\hat{\zeta}$ are harmonic in
$M\backslash\Gamma$ and continuously differentiable up to
$\Gamma$. In particular $u=u^{f}$ and
$\partial_{\nu}u=\Lambda_{g}f$, where $u^{f}$ is the solution of
(\ref{Eq 1}), (\ref{Eq 2}) and $\Lambda_{g}$ is the DN-map
of the above constructed $(M,g)$. Moreover, the Cauchy-Riemann
condition (\ref{CR--}) holds. Passing in (\ref{CR--}) to the trace
on $\Gamma$, one obtains
$$\Lambda_{g}f=\partial_{\nu}u=\partial_{\gamma}h=\Lambda f, \qquad \partial_{\nu}v=-\partial_{\gamma}f.$$
Since $f$ is arbitrary, we have proved that ${\rm
Ker}(\partial_{\gamma}+\Lambda J\Lambda)\subset{\rm
Ker}(\partial_{\gamma}+\Lambda_{g} J\Lambda_{g})$ and $\Lambda
f=\Lambda_{g}f$ for any $f\in {\rm Ker}(\partial_{\gamma}+\Lambda
J\Lambda)$.
\smallskip

\noindent${\bullet}$\,\,\,
Let us show that ${\rm Ker}(\partial_{\gamma}+\Lambda
J\Lambda)={\rm Ker}(\partial_{\gamma}+\Lambda_{g} J\Lambda_{g})$.
By {\bf iv}, the dimension $q$ of the factor-space ${\rm
Ker}(\partial_{\gamma}+\Lambda_{g} J\Lambda_{g})/{\rm
Ker}(\partial_{\gamma}+\Lambda J\Lambda)$ is finite. In view of
(\ref{Goperator}), one has
\begin{equation}
\label{kerups} {\rm Ker}\Upsilon=\{f+iJ\Lambda f+ic \ | \ f\in
{\rm Ker}(\partial_{\gamma}+\Lambda J\Lambda), \ c\in\mathbb{R}\}.
\end{equation}
Denote by $\mathfrak{A}^{\infty}$ the algebra of traces of all
holomorhic smooth functions on $M$; obviously, ${\rm Ker}\Upsilon$
is a subalgebra of $\mathfrak{A}^{\infty}$. From Cauchy-Riemann
conditions on $\Gamma$, the representation
$$\mathfrak{A}^{\infty}:=\{f+iJ\Lambda f+ic \ | \ f\in {\rm Ker}(\partial_{\gamma}+\Lambda_{g} J\Lambda_{g}), \ c\in\mathbb{R}\}.$$
is valid. Comparison of the last two formulas shows that
the algebra $\mathfrak{A}^{\infty}$ is a finite-dimensional
extension of the algebra ${\rm Ker}\Upsilon$, and dimension of the
factor-space $\mathfrak{A}^{\infty}/{\rm Ker}\Upsilon$ is equal to
$q$.

Suppose that $q>0$ and choose the elements
$\theta_{1},\dots,\theta_{q}\in \mathfrak{A}^{\infty}$ linearly
independent modulo ${\rm Ker}\Upsilon$. Then any
$\theta\in\mathfrak{A}^{\infty}$ can be represented as
\begin{equation}
\label{findimext rep}
\theta=\sum_{k=1}^{q}c_{k}(\theta)\theta_{k}+\tilde{\theta},
\end{equation}
where $c_{q}(\theta)\in\mathbb{C}$ and $\tilde{\theta}\in {\rm
Ker}\Upsilon$. Take any nonconstant $\eta\in{\rm Ker}\Upsilon$.
Representation (\ref{findimext rep}) implies
\begin{equation}
\label{findimext rep 1} \eta\theta_{l}=\sum_{k=1}^{q}
T_{kl}\,\theta_{k}+\tilde{\theta}_{l},
\end{equation}
where $T$ is a complex $q\times q-$matrix and
$\tilde{\theta}_{l}\in {\rm Ker}\Upsilon$. Choose an arbitrary
eigenpair $\lambda,X=(X_{1},\dots,X_{q})^{\,\rm tr}$ of $T$ and
denote
$$\Theta:=\sum_{k}X_{l}\theta_{l}, \quad \tilde{\Theta}:=\sum_{l=1}^{q}X_{l}\tilde{\theta}_{l}.$$
From (\ref{findimext rep 1}) it follows that
$$
\eta\Theta=\sum_{l=1}^{q}X_{l}\eta\theta_{l}=\sum_{k=1}^{q}\left(\sum_{l=1}^{q}T_{kl}X_{k}\right)\theta_{k}=
\lambda\sum_{k=1}^{q}X_{k}\theta_{k}+\sum_{l=1}^{q}X_{l}\tilde{\theta}_{l}=\lambda\Theta+\tilde{\Theta}.
 $$
Note that $\eta-\lambda e$ does not vanish identically on
any segment $\Gamma'$ of $\Gamma$ of non-zero length (indeed,
since $\hat{\eta}$ is holomorphic and smooth on $M$, the equality
$\eta=\lambda e$ on $\Gamma'$ implies $\eta=\lambda e$ on the
whole $\Gamma$). So,
\begin{equation}
\label{findimext rep 2} \Theta:=\frac{\tilde{\Theta}}{\eta-\lambda
e}
\end{equation}
holds on $\Gamma$, where both numerator and denominator are
elements of ${\rm Ker}\Upsilon$. Note that $X\ne 0$ and
$\Theta\notin{\rm Ker}\Upsilon$. Similarly, representation
(\ref{findimext rep}) yields
\begin{equation}
\label{findimext rep 3}
\Theta^{l}=\sum_{k=1}^{q}N_{kl}\theta_{k}+\tilde{\Theta}_{l},
\quad l=1,\dots, q,
\end{equation}
where $N$ is a complex $q\times q-$matrix and
$\tilde{\Theta}_{l}\in {\rm Ker}\Upsilon$.

If ${\rm det}N=0$, then there exists a non-zero
$Y=(Y_{1},\dots,Y_{q})^{\,\rm tr}\in {\rm Ker}N$ and the
polynomial $P(\Theta)=\sum_{l=1}^{q}Y_{l}\Theta^{l}$ admits the
representation
$$P(\Theta)=\sum_{k=1}^{q}\big(\sum_{l=1}^{q}N_{kl}Y_{l}\big)\theta_{k}+\sum_{l=1}^{q}Y_{l}\tilde{\Theta}_{l}=0+\sum_{l=1}^{q}Y_{l}\tilde{\Theta}_{l}.$$
Therefore $P(\Theta)\in {\rm Ker}\Upsilon$ and, due to
(\ref{findimext rep 2}) and condition {\bf ii}, one has $\Theta\in
{\rm Ker}\Upsilon$, which leads to a contradiction.

If ${\rm det}N\ne 0$ and $N'$ is the inverse matrix to $N$, then
(\ref{findimext rep 3}) implies
\begin{equation*}
\theta_{s}-\sum_{l=1}^{q}N'_{ls}\Theta^{l}=\sum_{l=1}^{q}N'_{ls}\tilde{\Theta}_{l}\in
{\rm Ker}\Upsilon.
\end{equation*}
This means that $\Theta,\Theta^{2}\dots,\Theta^{q}$ are linearly
independent modulo ${\rm Ker}\Upsilon$. So, we can assume that
$\theta_{k}=\Theta_{k}$. Now, formula (\ref{findimext rep})
renders
$$\mathcal{R}(\Theta):=\Theta^{q+1}-\sum_{k=1}^{q}c_{k}\theta_{k}\in {\rm Ker}\Upsilon,$$
where $c_{k}\in\mathbb{C}$. Since the polynomial $\mathcal{R}$ is
of degree $q+1>0$, the inclusion $\mathcal{R}(\Theta)\in {\rm
Ker}\Upsilon$, formula (\ref{findimext rep 2}) and condition {\bf
ii} yield  $\Theta\in {\rm Ker}\Upsilon$. This contradiction means
that $\mathfrak{A}^{\infty}={\rm Ker}\Upsilon$ and $q=0$. Thus, it
is proved that ${\rm Ker}\Upsilon$ is algebra of traces of all
holomorhic smooth functions on $M$ and ${\rm
Ker}(\partial_{\gamma}+\Lambda J\Lambda)={\rm
Ker}(\partial_{\gamma}+\Lambda_{g} J\Lambda_{g})$. In particular,
from (\ref{dimhandles}) it follows that ${\rm
dim\,}(\partial_{\gamma}+\Lambda
J\Lambda)C^{\infty}(\Gamma;\mathbb{R})=1-\mathscr{X}(M)$, where
$\mathscr{X}(M)$ is the Euler characteristic of $M$.
\smallskip

\noindent${\bullet}$\,\,\,
Thus, we have proved that $\Lambda$ coincides with DN-map
$\Lambda_{g}$ of the surface $(M,g)$ on the subspace
 \begin{equation}\label{Eq mathfrak K}
\mathfrak{K}:={\rm Ker}(\partial_{\gamma}+\Lambda J\Lambda)={\rm
Ker}(\partial_{\gamma}+\Lambda_{g} J\Lambda_{g})
 \end{equation}
of codimension $r:=1-\mathscr{X}(M)$ in
$C^{\infty}(\Gamma;\mathbb{R})$. To complete the proof of
sufficiency, it remains to show that $\Lambda
f_{1}=\Lambda_{g}f_{1},\dots,\Lambda f_{r}=\Lambda_{g}f_{r}$,
where $f_{1},\dots, f_{r}$ are some functions from
$C^{\infty}(\Gamma;\mathbb{R})$ linearly independent modulo
$\mathfrak{K}$. Before that, recall the terminology associated
with vector fields on the Riemannian manifolds and some well-known
facts.

The vector fields are the $TM_x$-valued functions on $M$ (the
cross-sections of $TM$). A field of the form $b=\nabla_g\varphi$
is called {\it potential}, $\varphi$ being a potential. A field
$a$ is {\it harmonic} if ${\rm div}_{g}\,a={\rm div}_{g}(\Phi
a)=0$ holds. The rotation $\Phi$ preserves the harmonicity. Each
harmonic field is {\it locally} potential. If $b=\nabla_g\varphi$
is harmonic then the potential $\varphi$ is a harmonic function:
$\Delta_g\varphi=0$, the opposite being also true.

So, let $f_{1},\dots, f_{r}$ be linearly independent modulo
$\mathfrak{K}$, $u_{j}$ the solution of problem (\ref{Eq 1}),
(\ref{Eq 2}) with $f=f_{j}$. The vector-fields
$a_{j}:=\Phi\nabla_g u_{j}$ are harmonic in $M$. Note that any
non-zero linear combination of $a_{j}$ is not a potential field in
$M$. Indeed, if $\sum_{j=1}^{r}c_{j}a_{j}=\nabla_g v$, then the
function $w:=u+iv$, where $u:=\sum_{j=1}^{r}c_{j}u_{j}$, is
holomorphic in $M$. Then $w|_{\Gamma}\in
\mathfrak{A}^{\infty}(M)={\rm Ker}\Upsilon$ and $\Re
w|_{\Gamma}=\sum_{j=1}^{r}c_{j}f_{j}\in \mathfrak{K}$ in view of
(\ref{kerups}). Since $f_{k}$ are linearly independent modulo
$\mathfrak{K}$, all $c_{j}$ equal zero.

Although $a_{j}$ are not potential on $M$, they can be represented
as gradients of some multi-valued functions $V_{j}$ which are
defined on an appropriate cover $\mathbb{M}$ of the surface $M$.
The cover $\mathbb{M}$ is constructed in the following way. Let
$\mathcal{D}$ be a surface diffeomorphic to an open disk in
$\mathbb{R}^{2}$ and $\partial{\mathcal D}=\Gamma$. Gluing up the
boundaries of $M$ and $\mathcal D$, we obtain the closed compact
surface $M'=M\cup\mathcal{D}$ of genus
$${\rm gen}\,M'=1-\frac{\mathscr{X}(M')}{2}=1-\frac{\mathscr{X}(M)+1}{2}=\frac{r}{2}.$$
As is well known, the metric tensor $g$ and rotation $\Phi$
on $M$ can be extended to the (smooth) metric tensor $g'$
and rotation $\Phi'$ on the whole $M'$.

Let $\mathbb{M}'$ be the universal covering of $M'$ (see, for
definition, \S 5, \cite{Forster}), which is a simply connected
Riemann surface, and $\pi': \mathbb{M}'\mapsto M'$ the projection,
which is a local homeomorphism. Tensor $g'$ and rotation $\Phi'$
on $M'$ induce the tensor ${\rm g}':=\pi'_{*}g'$ and the rotation
$\tilde\Phi:=\pi'_{*}\Phi'$ on $\mathbb M'$. As a result, $\pi':
(\mathbb{M}',{\rm g}')\mapsto (M',g')$ turns out to be a local
isometry. At last, we get the required cover for $(M,g,\Phi)$ as
the collection $(\mathbb{M},\pi,{\rm g},{\dot\Phi})$, where
$\mathbb{M}:=\mathbb{M}'\backslash{\pi'^{\,\,-1}(\mathcal{D})}$ is
the surface with the boundary
$\partial\mathbb{M}=\pi'^{\,\,-1}(\Gamma)$, ${\rm g}:={\rm
g'}|_{\mathbb{M}}$, ${\dot\Phi}:=\tilde\Phi|_{\mathbb{M}}$, and
$\pi:=\pi'|_{\mathbb{M}}$.

Recall that the solutions $u_j$ and fields $a_j$ correspond to the
functions $f_{1},\dots, f_{r}$ which are linearly independent
modulo $\mathfrak{K}$ (see (\ref{Eq mathfrak K})). Introduce the
vector fields $A_{j}:=\pi_{*}a_{j}=\dot\Phi\nabla_{\rm g}
(u_{j}\circ\pi)$ on $\mathbb{M}$ and the functions
$$\mathbb{M}\ni x\mapsto V_{j}(x)=\int_{\mathcal{L}}{\rm g}(A_{j},l)\,dl\in\mathbb{R},$$
where $\mathcal{L}$ is an arbitrary curve in $\mathbb{M}$ that
connects a fixed point $x_{0}\in\mathbb{M}$ to a point $x$. In
what follows, we denote by $l$ and $dl$ the unit tangent vector
and the length element on the curve, respectively. Since
$\mathbb{M}=\mathbb{M}'\backslash\pi'^{\,\,-1}(\mathcal{D})$ is no
longer simply connected, one needs to check that $V_{j}$ are
single-valued on $\mathbb{M}$. To this end, it suffices to show
that $\int_{\tilde{\Gamma}}{\rm g}(A_{j},l)\,dl=0$ for any
connected component $\tilde{\Gamma}$ of $\pi^{\,\,-1}(\Gamma)$.
Since $\tilde{\Gamma}$ is isometric to $\Gamma$, one needs to
check only that $\int_{\Gamma}g(a_{j},\gamma)\,d\gamma=0$. By the
Green formula, one has
$$\int_{\Gamma}g(a_{j},\gamma)\,d\gamma=\int_{\Gamma}g(\Phi\nabla_g u_{j},\gamma)\,d\gamma=\int_{\Gamma}\partial_{\nu}u_{j}\,d\gamma=\int_{M}\Delta_g u_{j}\,dx=0$$
in view of harmonicity of $u_j$. So, we have provided the
functions $V_{j}$ such that $\nabla_{\rm g}
V_{j}=A_{j}=\dot\Phi\nabla_{\rm g} (u_{j}\circ\pi)$ holds on
$\mathbb{M}$. This means that the functions
 \begin{equation}\label{Eq W_j}
W_{j}:=u_{j}\circ\pi+i\,V_{j},\qquad j=1,\dots, r
 \end{equation}
are holomorphic on $(\mathbb{M},{\rm g})$, whereas the
Cauchy-Riemann conditions $\nabla_{\rm g}\Im
W_j=\dot\Phi\nabla_{\rm g}\Re W_j$ hold.
\smallskip

\noindent${\bullet}$\,\,\,
We are going to show that the functions $f_{1},\dots, f_{r}$ can
be chosen in such a way that $e^{W_{j}}=w_{j}\circ\pi$, where
$w_{j}$ are holomorphic functions in $M$.

Introduce the groups
\begin{align*}
& {\rm Deck}(\mathbb{M}/M)\,\,\, :=\{\phi \ | \ \phi \text{ is automorphism of } \mathbb{M}, \ \pi\circ\phi=\pi\},\\
& {\rm Deck}(\mathbb{M}'/M') :=\{\phi' \ | \ \phi \text{ is
automorphism of } \mathbb{M}', \ \pi'\circ\phi'=\pi'\}
\end{align*}
of fiber-wise automorphisms of $\mathbb{M}$ and $\mathbb{M'}$,
respectively (see, e.g., \cite{Forster}, 5.4). Obviously, if
$\phi'\in {\rm Deck}(\mathbb{M}'/M')$,
then 
$\phi'|_\mathbb M\in{\rm Deck}(\mathbb{M}/M)$. Conversely, if
$\phi\in {\rm Deck}(\mathbb{M}/M)$, then it can be lifted to
$\phi'\in{\rm Deck}(\mathbb{M}'/M')$ such that
$\phi'|_{\mathbb{M}}=\phi$. Indeed, if $x$ belongs to a connected
component $\tilde{\mathcal{D}}$ of
$\pi'^{\,\,-1}(\mathcal{D})$, then $\phi'(x)$ is uniquely
determined by its projection $\pi'(\phi'(x))=\pi'(x)$ and by the
fact that the boundary of the connected component of
$\pi'^{\,\,-1}(\mathcal{D})$ containing $\phi'(x)$ must coincide
with $\phi(\partial\tilde{\mathcal{D}})$. So, the map
$\phi'\mapsto\beta\phi'=\phi'|_{\mathbb{M}}$ is an isomorphism of
the groups ${\rm Deck}(\mathbb{M}'/M')$ and ${\rm
Deck}(\mathbb{M}/M)$.

Denote by $\pi_{1}(M')$ the fundamental group of $M'$ and
by $[L]$ the homotopy class of a closed curve $L$ in $M'$. In view of
Proposition 5.6, \cite{Forster}, the groups $\pi_{1}(M')$ and ${\rm Deck}(\mathbb{M}'/M')$ are
isomorphic. The isomorphism
$$\alpha: \ {\rm Deck}(\mathbb{M}'/M')\mapsto
\pi_{1}(M')$$ is constructed as follows. Let $\phi'\in {\rm
Deck}(\mathbb{M}'/M')$. Choose an arbitrary point $x\in
\mathbb{M'}$ and a curve $\mathcal{L}_{\phi'}$ which connects $x$
to $\phi'(x)$. Then $\pi'(\mathcal{L}_{\phi'})$ is a closed curve
in $M'$ due to the equality $\pi'(\phi'(x))=\pi'(x)$. It turns out
that the homotopy class $[\pi'(\mathcal{L}_{\phi'})]$ of the curve
$\pi'(\mathcal{L}_{\phi'})$ does not depend on the choice of
$\mathcal{L}_{\phi'}$ and $x$. The required isomorphism $\alpha$
is defined by the rule
$$\alpha(\phi'):=[\pi'(\mathcal{L}_{\phi'})].$$
Also, $\alpha\circ\beta^{-1}$ is an isomorphism of
groups ${\rm Deck}(\mathbb{M}/M)$ and $\pi_{1}(M')$.

Since $M'$ is a surface of the genus ${\rm gen}M'=r/2$, there are
$2\,{\rm gen}M'=r$ generators $[L_{1}],\dots,[L_{r}]$ of the
fundamental group $\pi_{1}(M')$. Note that, since $\mathcal{D}$ is
simply connected, one can deform the curves $L_{j}$, preserving
their homotopy class, in such a way that any $L_{j}$ does not
intersect $\mathcal{D}$. Thus, we assume that
$L_{1},\dots,L_{r}\subset M$. Since the groups ${\rm
Deck}(\mathbb{M}/M)$ and $\pi_{1}(M')$ are isomorphic, the
automorphisms $\phi_{j}:=\beta\circ\alpha^{-1}([L_{j}])$,
$j=1,\dots,r$ generate the group ${\rm Deck}(\mathbb{M}/M)$.
Therefore, a function $V$ on $\mathbb M$ can be represented as
$V=v\circ\pi$ if and only if $V\circ\phi_{j}=V$, $j=1,\dots,r$.

Suppose that $V$ is a function on $M$ such that $\nabla_{\rm g}
V=A:=\pi_{*}a$, where $a$ is a vector field on $M$. Then
\begin{equation*}
V(\phi_{j}(x))-V(x)=\int_{\mathcal{L}_{j}}{\rm g}(A,l)\,dl,
\end{equation*}
where $\mathcal{L}_{j}$ connects $x$ to $\phi_{j}(x)$. Since the
field $A=\pi_{*}a$ is invariant under action of the group ${\rm
Deck}(\mathbb{M}/M)$, the right-hand side does not depend on $x$
and one can choose $\mathcal{L}_{j}$ to provide
$\pi(\mathcal{L}_{j})=L_{j}$. Then the difference
$V(\phi_{j}(x))-V(x)$ is equal to
$$T_{j}(a):=\int_{L_{j}}{\rm g}(a,l)\,dl.$$
Thus, $V=v\circ\pi$ and $a=\nabla_g v$ if and only if
$T_{1}(a)=\dots=T_{r}(a)=0$.

Introduce the $r\times r$-matrix $T$ with entries
$T_{ij}=T_{i}(a_{j})$. Recall that any non-zero linear combination
$\sum_{j=1}^{r}c_{j}a_{j}$ is not potential field in $M$. This
means that all
$T_{i}(\sum_{j=1}^{r}c_{j}a_{j})=\sum_{j}T_{ij}c_{j}$ are zero if
and only if $c_{1}=\dots=c_{r}=0$. Thus, $T$ is invertible. Denote
$f'_{s}:=2\pi\sum_{l=1}^{r}R_{ls}f_{l}$, where $R=T^{-1}$. Then
$f'_{1},\dots,f'_{r}$ are linear independent modulo
$\mathfrak{K}$. Introduce the new functions
$$V'_{s}=2\pi\sum_{l=1}^{r}R_{ls}V_{l}, \quad W_{s}=2\pi\sum_{l=1}^{r}R_{ls}W_{l},$$
that are determined by $f'_{s}$ in the same way as $V_{s}$ and
$W_{s}$ are determined by $f_{s}$ (see (\ref{Eq W_j})). Then
$\nabla_{\rm g}
V'_{s}=2\pi\sum_{l=1}^{r}R_{ls}A_{l}=\pi_{*}a'_{s}$, where
$a'_{s}=2\pi\sum_{l=1}^{r}R_{ls}a_{l}$ and
$$T_{j}(a'_{s})=2\pi\sum_{l=1}^{r}T_{jl}R_{ls}=2\pi\delta_{js}$$
holds. By the latter, one has
$$V'_{s}\circ\phi_{j}-V'_{s}=W'_{s}\circ\phi_{j}-W'_{s}=2\pi\delta_{js}, \qquad  j=1,\dots,r.$$
Hence,
$$e^{W'_{s}\circ\phi_{j}}=e^{W'_{s}}$$
for any $j,s=1,\dots,r$. This means that $e^{W'_{s}}$ can be
represented as $e^{W'_{s}}=w_{s}\circ\pi$, where $w_{s}$ is a
function on $M$. Since $W'_{s}$ is holomorphic in $\mathbb{M}$,
the function $w_{s}$ is holomorphic in $M$. Replacing $f_{s}$ by
$f'_{s}$ (what is the same, omitting `prime' everywhere in the
notation), one obtains $e^{W_{j}}=w_{j}\circ\pi$.

So, we have provided the functions $f_{1},\dots,f_{r}$ with the
properties claimed at the beginning of the paragraph.
\smallskip

\noindent${\bullet}$\,\,\,
Since $w_j$ is holomorphic on $M$, one has
$\zeta_{j}:=w_{j}|_{\Gamma}\in{\rm Ker\,}\Upsilon$. Also, one has
 \begin{align*}
{\rm log}|w_{j}(\pi(x))|=\Re W_{j}(x)=u_{j}(\pi(x)), \qquad
x\in\mathbb{M}.
 \end{align*}
In particular,
$${\rm log\,}|w_{j}|=u_{j}\quad\text{and}\quad{\rm log\,}|\zeta_{j}|=f_{j}$$
holds on $M$ and $\Gamma$ respectively.

Since $W_{j}$ is holomorphic on $\mathbb{M}$, the Cauchy-Riemann
conditions yield
 \begin{align*}
&
(\partial_{\nu}u_{j})\circ\pi=\partial_{\nu}(u_{j}\circ\pi)=\partial_{\nu}\,\Re
W_{j}=
\partial_{\gamma}\Im W_{j}=\partial_{\gamma}\Im\,{\rm
log\,}(w_{j}\circ\pi)=\\
& =
\partial_{\gamma}{\rm arg\,}(w_{j}\circ\pi)=(\partial_{\gamma}{\rm arg\,}w_{j})\circ\pi
 \end{align*}
on $\pi^{-1}(\Gamma)$. This means that
\begin{equation}
\label{end 1}
\Lambda_{g}f_{j}=\partial_{\nu}u_{j}=\partial_{\gamma}{\rm
arg\,}\zeta_{j}
\end{equation}
holds on $\Gamma$. In the mean time, by virtue of {\bf vii} we
have
\begin{equation}
\label{end 2} \Lambda f_{j}=\Lambda\, {\rm
log\,}|\zeta_{j}|=\partial_{\gamma}{\rm arg\,}\zeta_{j}.
\end{equation}
Comparing (\ref{end 1}) and (\ref{end 2}), one obtains $\Lambda
f_{j}=\Lambda_{g}f_{j}$ for any $j=1,\dots,r$. Together with what
was proved above, this means that $\Lambda=\Lambda_{g}$ and,
hence, $\Lambda$ is the DN-map of the surface $(M,g)$.
\smallskip

{\it The sufficiency of the conditions {\bf i}-{\bf vii} is
established.}
\medskip

Theorem \ref{main} is proved.

\end{document}